\documentclass[12pt,twoside,reqno]{amsart}
\usepackage[colorlinks=true,citecolor=blue]{hyperref}
\usepackage{mathptmx, amsmath, amssymb, amsfonts, amsthm, mathptmx, enumerate, color,cases,graphics}
\usepackage{graphicx}
\usepackage{subfig}
\usepackage{epstopdf}
\usepackage{multirow}
\usepackage{lineno}
\usepackage{cite}
\usepackage[mathcal]{eucal}
\usepackage{tabularx}
\newtheorem{theorem}{Theorem}[section]

\newtheorem{lemma}{Lemma}[section]
\newtheorem{proposition}{Proposition}[section]
\theoremstyle{definition}

\newtheorem{remark}{Remark}[section]

\newtheorem{assumption}{Assumption}
\newtheorem{algorithm}{Algorithm}[section]
\numberwithin{equation}{section}

\begin{document}
\setcounter{page}{1}

\vspace*{1.0cm}
\title[modified PRP-type conjugate gradient method]
{A modified Polak-Ribi\`{e}re-Polyak type conjugate gradient method with two stepsize strategies for vector optimization}
\author[Y.S. Bai,J.W. Chen]{Yushan Bai$^{1}$,Jiawei Chen$^{1}$,Kaiping liu$^{1}$}
\maketitle
\vspace*{-0.6cm}

\begin{center}
{\footnotesize {\it

$^1$School of Mathematics and statistics, Southwest University,  Chongqing 400715, China\\

}}\end{center}
\vskip 4mm

{\small\noindent {\bf Abstract.}

In this paper, in order to find critical points of vector-valued functions with respect to the partial order induced by a closed, convex, and pointed cone with nonempty interior, we propose a nonlinear modified Polak-Ribi\`{e}re-Polyak type conjugate gradient method with a nonnegative conjugate parameter. We show that the search direction in our method satisfies the sufficient descent condition independent of any line search. Furthermore, under mild assumptions, we obtain the results of global convergence with the standard Wolfe line search conditions as well as the standard Armijo line search strategy without convexity assumption of the objective functions. Computational experiments are given to show the effectiveness of the proposed method.

\noindent {\bf Keywords.}
Vector optimization; Conjugate gradient direction; Line search strategy;  Pareto critical point}

\renewcommand{\thefootnote}{}
\footnotetext{
E-mail addresses: YushanBai3@outlook.com, j.w.chen713@163.com, kaipliu@163.com.

\par
}
\section{Introduction}

Let us first consider the single-objective problem
\begin{equation}
\text{minimize} \quad f(x),\quad x\in\mathbb{R}^n,
\end{equation}
where $f:\mathbb{R}^n\rightarrow \mathbb{R}$ is continuously differentiable. The most classic first-order method to solve this problem is the steepest descent method, while the most widely used second-order method is the Newton's method. And the conjugate gradient method proposed in \cite{3} is one of the most commonly used and effective optimization methods between the steepest descent method and the Newton's method, which has the characteristics of fast convergence rate, lower memory requirement and less computation cost, and has important applications in both linear and nonlinear optimization. Formally, the conjugate gradient method generates a sequence $\{x^k\}$ given by
\begin{equation}
x^{k+1}=x^k+\alpha_k d^k,\quad k\geq 0,
\end{equation}
where the stepsize $\alpha_k>0$ is obtained by a line search strategy and the search direction $d^k$ is defined by
\begin{equation}\label{sdk}
d^{k}=\begin{cases}
-\nabla f(x^k), & \text { if } k=0, \\
-\nabla f(x^k)+ \beta_{k} d^{k-1}, & \text { if } k \geq 1,
\end{cases}
\end{equation}
where $\beta_k$ is a scalar algorithmic parameter. For nonquadratic functions, Different choices for the conjugate parameter $\beta_k$ in \eqref{sdk} result in different algorithms, known as \emph{nonlinear conjugate gradient methods}. Although too many to name, some notable choices would include:
\begin{equation*}
\begin{aligned}
&\text{Fletcher-Reeves(FR) \cite{3}}:\quad \beta_k=\frac{\langle g^k,g^k\rangle}{\langle g^{k-1},g^{k-1}\rangle};\\
&\text{Conjugate descent(CD) \cite{4}}:\quad \beta_k=-\frac{\langle g^k,g^k\rangle}{\langle d^{k-1},g^{k-1}\rangle};\\
&\text{Dai-Yuan(DY) \cite{5}}:\quad \beta_k=\frac{\langle g^k,g^k\rangle}{\langle d^{k-1},g^k-g^{k-1}\rangle};\\
&\text{Polak-Ribi\`{e}re-Polyak(PRP) \cite{6}}:\quad \beta_k=\frac{\langle g^k,g^k-g^{k-1}\rangle}{\langle g^{k-1},g^{k-1}\rangle};\\
&\text{Hestenes-Stiefel(HS) \cite{7}}:\quad \beta_k=\frac{\langle g^k,g^k-g^{k-1}\rangle}{\langle d^{k-1},g^k-g^{k-1}\rangle},
\end{aligned}
\end{equation*}
where $g^k=\nabla f(x^k)$ and $\langle \cdot,\cdot \rangle$ denotes the usual inner product. We expect to get conjugate parameters $\beta_k$ that make $d^k$ descent directions in the sense of $\langle g^k,d^k\rangle<0$ for all $k\geq 0$, or make $d^k$ meet a more stringent condition, which is called sufficient descent condition and defined as
\begin{equation}\label{ssdc}
\langle g^k,d^k\rangle\leq-c\Vert g^k\Vert^2,
\end{equation}
for some $c>0$ and all $k\geq 0$, where $\Vert \cdot\Vert$ denotes the Euclidian norm. A significant advantage of the FR, CD, and DY methods is that if a line search satisfying the Wolfe conditions is used, the corresponding search directions are verified to be descent. However, the PRP and HS methods do not necessarily generate descent directions even when Wolfe line searches are employed.

In this paper, we consider the following unconstrained vector optimization problem(VOP)
\begin{equation}\label{1}
\text{minimize}_{K} \quad F(x), \quad x\in\mathbb{R}^n,
\end{equation}
where $F:\mathbb{R}^{n}\to \mathbb{R}^{m}$ is continuously differentiable, and $K\subset \mathbb{R}^m$ is a pointed, closed and convex cone with nonempty interior $\text{int} (K)$. The partial order $\preceq_K(\prec_K)$ is given by $u\preceq_K v(u \prec_K v)$ if and only if $v-u \in K (v-u\in \text{int}(K))$. In vector optimization, the concept of optimality is replaced by the concept of Pareto optimality or efficiency, and we seek to find \emph{K-Pareto optimal point} or \emph{K-efficient point}. In practical applications, we usually take $K=\mathbb{R}_+^m$, then Problem \eqref{1} corresponds to the multiobjective optimization problem.

Vector optimization problems are a significant extension of multiobjective optimization, and there are a large number of real life applications of multicriteria and vectorial optimization, such as engineering design \cite{8,11}, finance \cite{9},  machine learning \cite{10}, space exploration \cite{12}, management science \cite{13,14}, environmental analysis and so on. Due to the wide application of vector-valued optimization, the development of strategies for solving vector-valued optimization problems has attracted wide attention. At present, a lot of research has been made in the theory and algorithm of solving vector optimization problems, and the common methods for solving vector optimization problems include scalarization approaches and descent methods. The scalarization approaches for solving vector optimization problems are to convert the original vector optimization problems into the parameterized single objective ones; see \cite{48,49}. The drawback of this method is that even when the original vector-valued problem has solutions, the selection of parameters may lead to unbounded numerical problems (and thus unsolvable). However, the descent methods do not require any parameter information, and thus usually perform better in numerical experiments.  Many descent methods for solving scalar optimization problems have been extended to vector-valued optimization, such as projected gradient method \cite{15,16,17}, Newton method \cite{18,19,20}, steepest descent method \cite{1,22}, proximal point method \cite{23,24} and so on.

In recent years, the conjugate gradient method has been extended from solving single-objective problems to solving multiobjective \cite{25} and vector-valued problems \cite{2,26,27}. The first work in this line was \cite{2}, Lucambio P\`{e}rez and Prudente proposed a nonlinear conjugate gradient algorithm (NLCG), and generate a sequence of iterates by the following form:
\begin{equation}
x^{k+1}=x^k+\alpha_k d^k,\quad k\geq 0,
\end{equation}
where the stepsize $\alpha_k>0$ is obtained by standard Wolfe or strong Wolfe line search strategy and the search direction $d^k$ is defined by
\begin{equation}
d^{k}=\begin{cases}
v(x^k), & \text { if } k=0, \\
v(x^k)+ \beta_{k} d^{k-1}, & \text { if } k \geq 1,
\end{cases}
\end{equation}
where $\beta_k$ is scalar algorithmic parameters, which are extended from the FR, CD, DY, PRP, and HS conjugate gradient algorithms for the single-objective case, and they are defined as follows
\begin{equation*}
\begin{aligned}
&\text{Fletcher-Reeves(FR)}:\quad \beta_k=\frac{h(x^k,v(x^k))}{h(x^{k-1},v(x^{k-1}))};\\
&\text{Conjugate descent(CD)}:\quad \beta_k=\frac{h(x^k,v(x^k))}{h(x^{k-1},d^{k-1})};\\
&\text{Dai-Yuan(DY)}:\quad \beta_k=-\frac{h(x^k,v(x^k))}{h(x^{k},d^{k-1})-h(x^{k-1},d^{k-1})};\\
&\text{Polak-Ribi\`{e}re-Polyak(PRP)}:\quad \beta_k=\frac{-h(x^k,v(x^k))+h(x^{k-1},v(x^k))}{-h(x^{k-1},v(x^{k-1}))};\\
&\text{Hestenes-Stiefel(HS)}:\quad \beta_k=\frac{-h(x^k,v(x^k))+h(x^{k-1},v(x^k))}{h(x^{k},d^{k-1})-h(x^{k-1},d^{k-1})},
\end{aligned}
\end{equation*}
where $h(\cdot,\cdot)$ is defined in the next section. In \cite[Theorem 5.11]{2}, by assuming that the search direction $d^k$ is a $K$-descent direction, Lucambio P\`{e}rez and Prudente established the convergence result related to the PRP+ parameter given by $\beta_k^{PRP+}:=\max\{\beta_k^{PRP},0\}$, which remind us that the nonegativeness of parameter $\beta_k$ seems to be essential for obtaining the convergence result of the conjugate gradient method, while the nonegativeness of the PRP parameter $\beta_k^{PRP}$ cannot be guaranteed. To address this weakness of the PRP parameter, we extend the method considered in \cite{28} to the vector context because the parameter $\beta_k$ of this method are nonnegative and this method show superior performance in numerical experiments. We propose a nonlinear modified Polak-Ribi\`{e}re-Polyak type conjugate gradient method with a nonnegative conjugate parameter, which is extended from \cite{28}, and we show that the search direction in our proposed method satisfies the sufficient descent condition no matter what stepsize strategy is adopted. Furthermore, under mild assumptions, which are natural extensions of those made for the scalar case, we obtain the results of global convergence with the standard Wolfe line search as well as the standard Armijo line search without convexity assumption of the objective functions.

The paper is organized as follows. In the next section, we present some notations, definitions and preliminary results. In Section 3, we propose the modified PRP-type conjugate gradient method and investigate some properties of this method. The global convergence of the full sequence generated by the proposed method with the standard Wolfe line search as well as the standard Armijo line search is provided in Section 4. In section 5, some numerical experiments are reported to show the ability of the proposed method. Finally, we give some concluding remarks in Section 6.

\section{Preliminaries}

In this section, we present the vector optimization problem studied in the present work, the first order optimality condition for it, and some notations and defnitions.  Throughout this paper, let $\langle\cdot,\cdot\rangle$ stands for the inner product in $\mathbb{R}^n$ and $\Vert \cdot\Vert$ denotes the norm, that is $\Vert x\Vert=\sqrt{\langle x,x\rangle }$ for $x\in\mathbb{R}^n$. And we denote the convex hull of $A\subset\mathbb{R}^m$ by $\text{conv}(A)$,  and the cone of $A$ by $\text{cone}(A)$, let $K\subset \mathbb{R}^m$ be a pointed, closed and convex cone, with nonempty interior $\text{int} (K)$. The partial induced by $K$, $\preceq_K$ is defined as follows
\begin{equation*}
u \preceq_K v \quad\text{if and only if} \quad v-u\in K,
\end{equation*}
and the partial induced by $\text{int}(K)$, $\prec_K$ is defined as follows
\begin{equation*}
u \prec_K v \quad\text{if and only if} \quad v-u\in \text{int}(K).
\end{equation*}

A point $x^{*}\in \mathbb{R}^n$  is called a \emph{K-Pareto optimal point} (or \emph{K-Pareto point}) of \eqref{1} on $\mathbb{R}^n$, if there exists no other point $x\in \mathbb{R}^n$, such that $F(x)\preceq_K F(x^{*})$ and $F(x)\neq F(x^{*})$. The set of the objective values of all Pareto optimal solutions is also called Pareto frontier. In turn, a point $x^{*}\in \mathbb{R}^n$  is called a \emph{K-weak Pareto optimal point} (or \emph{K-weak Pareto}) of \eqref{1} on $\mathbb{R}^n$, if there exists no other point $x\in \mathbb{R}^n$, such that $F(x)\prec_K F(x^{*})$. It is clear that a Pareto optimal point is also a weak Pareto optimal point but not vice versa.
Since $F$  is continuously differentiable, the subdifferential of $F$ at $x\in \mathbb{R}^n$  coincides with the Jacobian of $F$, and is denoted by $J F(x)$ and the image of the Jacobian of $F$ at a point $x$ is denoted by $\text{Im}(J F(x))$. A first order optimality condition (necessary but in general not sufficient) for the problem \eqref{1} of a point $x\in \mathbb{R}^n$ is given by
\begin{equation}\label{oc}
(-\text{int}(K)) \cap \operatorname{Im}(J F(x))=\emptyset.
\end{equation}
which means that, for any $d\in \mathbb{R}^n$, we have $J F(x)d\notin -\text{int}(K)$. A point $x\in \mathbb{R}^n$  satisfying \eqref{oc} is called a \emph{K-Pareto critical point} or a \emph{K-stationary point} of problem \eqref{1}. Note that if $x\in \mathbb{R}^n$ is not a K-Pareto critical point, then there exists a direction $d\in \mathbb{R}^n$ satisfying $J F(x)d\in -\text{int}(K)$. This implies that $d$ is a $\emph{K-descent direction}$ for $F$ at $x$, i.e., there exists $\epsilon > 0$ such that $F(x+\alpha d)\prec_K F(x)$ for all $\alpha \in ]0,\epsilon]$.

The positive polar cone of $K\subset \mathbb{R}^m$ is the set
\begin{equation}
K^*:=\{\omega\in\mathbb{R}^m\mid\langle y,w\rangle\geq0,\forall y\in K\},
\end{equation}
since the set $K$ is closed and convex, $K=K^{**}$, and thus $-K=\{y\in\mathbb{R}^m\mid\langle y,w\rangle\leq0,\forall \omega\in K^*\}$ and $-\text{int}(K)=\{y\in\mathbb{R}^m\mid\langle y,w\rangle<0,\forall \omega\in K^*\backslash\{0\}\}$.
According to \cite[Remark 1.6]{44}, every cone in finite dimensional spaces has a closed convex bounded base if and only if it is pointed closed, which means that there is a compact set $C\subset\mathbb{R}^m$ satisfying
\begin{equation}\label{C_condition}
0\notin C\quad \text{and} \quad K^*=\text{cone}(\text{conv}(C)),
\end{equation}
since $\text{int}(K)\neq \emptyset$ and $C\subset K^*\backslash\{0\}$, it follows that $0\notin \text{conv}(C)$. Then
\begin{equation}\label{-k}
-K=\{y\in\mathbb{R}^m\mid\langle y,\omega\rangle\leq 0,\quad\forall \omega\in C\},
\end{equation}
and
\begin{equation}
-\text{int}(K)=\{y\in\mathbb{R}^m\mid\langle y,\omega\rangle<0,\quad\forall \omega\in C\}.
\end{equation}
\begin{remark}
It is known that $K=R_+$, and $C=\{1\}$ in single-objective optimization. As for multiobjective optimization, $K$ and $K^*$ are the positive orthant of $\mathbb{R}^m$ and we may take $C$ as the canonical basis of $\mathbb{R}^m$. If $K$ is a polyhedral cone, $C$ may be taken as a finite set of extremal rays of $K^*$.
\end{remark}
For generic $K$ which is a pointed, closed and convex cone with nonempty interior, we define
\begin{equation}
C:=\{w\in K^*\mid\Vert w\Vert=1\},
\end{equation}
then $C$ satisfies the condition that we mentioned in \eqref{C_condition}.

For the convenience of the subsequent description, we define $\phi:\mathbb{R}^m\to \mathbb{R}$ as follows
\begin{equation}\label{phiy}
\phi(y)=\sup\{\langle y,w \rangle \mid w\in C \},
\end{equation}
consider the compactness of $C$, the function $\phi$ is well defined. Then $-K$ and $-\text{int}(K)$ can be rewrite as follows
\begin{equation}
-K=\{y\in\mathbb{R}^m\mid\phi(y)\leq 0\},
\end{equation}
and
\begin{equation}
-\text{int}(K)=\{y\in\mathbb{R}^m\mid\phi(y)< 0\}.
\end{equation}

In the following Lemma, we will give some basic properties of the function $\phi$ stated in \cite[Lemma 3.1]{1}.
\begin{lemma}
From the definition of $\phi$, for $\forall$ $y,y^{\prime}\in \mathbb{R}^m$, the following statements hold:
\begin{enumerate}
\item[\rm(i)] $\phi(y+y^{\prime})\leq\phi(y)+\phi(y^{\prime})$ and $\phi(y)-\phi(y^{\prime})\leq\phi(y-y^{\prime})$;
\item[\rm(ii)] If $y\preceq_K y^{\prime}$, then $\phi(y)\leq \phi(y^{\prime})$; if $y\prec_K y^{\prime}$, then $\phi(y)< \phi(y^{\prime})$;
\item[\rm(iii)] $\phi$ is Lipschitz continuous with constant $1$.
\end{enumerate}
\end{lemma}
Note that $y\in K$ implies that $\phi(y)\geq 0$, and $y\in \text{int}(K)$ implies that $\phi(y)> 0$.

For the convenience of the subsequent description, we define $h:\mathbb{R}^n\times \mathbb{R}^n\to \mathbb{R}$ as follows
\begin{equation}\label{fxd}
h(x,d)=\phi(JF(x)d)=\sup \{ \langle JF(x)d,w \rangle\mid w\in C\} .
\end{equation}
From the definition of $h$, we know that $h$  can express $K$-Pareto critical point and $K$-descent direction of vector optimization problem, and we state it in the following Lemma.
\begin{lemma}\cite{1}
Let $d\in \mathbb{R}^n$, $x\in \mathbb{R}^n$ and $h(x,d)$ be defined as \eqref{fxd}, we have
\begin{enumerate}
\item[\rm(i)] $d$ is a $K$-descent direction at $x$ iff $h(x,d)<0$ ;

\item[\rm(ii)]$x$ is a $K$-Pareto critical point iff $h(x,d)\geq 0$ for any $d$.
\end{enumerate}
\end{lemma}

Drummond and Svaiter \cite{1} defined the steepest descent direction for vector optimization problem using the unique optimal solution of the following problem as
\begin{equation}\label{sd}
\min_{d\in\mathbb{R}^n} \quad h(x,d)+\frac{1}{2}\Vert d\Vert^{2}.
\end{equation}
Since $h(x,d)$ is a real closed convex function, the solution for \eqref{sd} exists and unique, we assume that $v(x)$ and $\theta(x)$ are the optimal solution and the optimal value of Problem \eqref{sd} from now on, for each $x\in \mathbb{R}^n$ respectively. That is
\begin{equation}\label{vk}
v(x)=\mathop{\arg\min}\limits_{d\in\mathbb{R}^n} \{h(x,d)+\frac{1}{2}\Vert d\Vert^{2}\},
\end{equation}
and
\begin{equation}\label{thetak}
\theta(x)=h(x,v(x))+\frac{1}{2}\Vert v(x)\Vert^{2}.
\end{equation}

Let us now state some basic results relating to the stationarity of a given point $x$ about $v(x)$ and $\theta(x)$.
\begin{lemma}\cite[Lemma 3.3]{1}\label{stop}
Consider $v(x)$ and $\theta(x)$ be defined as \eqref{vk} and \eqref{thetak}, The following statements hold true:
\begin{enumerate}
 \item[\rm(1)] For each $x\in\mathbb{R}^n$, $\theta(x)\leq0$;
 \item[\rm(2)] The following conditions are equivalent:
 \begin{enumerate}
   \item[\rm(i)] $x$ is not a $K$-Pareto critical point;
   \item[\rm(ii)] $v(x)\neq 0$;
   \item[\rm(iii)] $\theta(x)<0$.
 \end{enumerate}
 \item[\rm(3)] $v(\cdot)$ and $\theta(\cdot)$ are continuous.
\end{enumerate}
\end{lemma}
\begin{remark}
If $x\in \mathbb{R}^n$ is not a $K$-Pareto critical point, then we have
\begin{equation*}
h(x,v(x))\leq - \frac{\Vert v(x)\Vert^2}{2}<0,
\end{equation*}
and $v(x)$ is a $K$-descent direction for $F$ at $x$.
\end{remark}

\begin{lemma}\label{lema+b}\cite[Lemma 2.4]{2}
For any scalars $a,b$ and $\xi\neq0$, we have
\begin{equation}
(a+b)^2\leq (1+2\xi^2)a^2+(1+\frac{1}{2\xi^2})b^2.
\end{equation}
\end{lemma}

\section{A modified PRP-type conjugate gradient method and its property}
In this section, we will describe the modified PRP-type conjugate gradient method for vector optimization and then present some results which shows that our algorithm is well-defined.
\begin{algorithm}\label{alg1}
{\bf [modified PRP-type conjugate gradient method]}
\begin{enumerate}
 \item[{\rm Step 1.}] Let $x^{0}\in \mathbb{R}^n$ be an arbitrary initial point. Choose parameters $0<\rho<\sigma<1$, $\mu>2$, set $k=0$.

 \item[{\rm Step 2.}]Compute the direction $v(x^k)=\mathop{\arg\min}\limits_{d\in \mathbb{R}^n} \varphi_k(d)$, where $\varphi_k(d)=h(x^{k},d)+\frac{1}{2}\Vert d\Vert^{2}$.

 \item [{\rm Step 3.}] If $v(x^{k})=0$, STOP.  Otherwise, proceed to Step 4.

 \item [{\rm Step 4.}] Computing
\begin{equation}\label{dk}
d^{k}=\begin{cases}
v(x^{k}), & \text { if } k=0, \\
v(x^{k})+ \beta_{k}^{MPRP} d^{k-1}, & \text { if } k \geq 1,
\end{cases}
\end{equation}
where
\begin{equation}\label{betak}
\beta_{k}^{MPRP}=\frac{-h(x^k,v(x^k))(\vert h(x^{k-1},v(x^k))\vert+h(x^{k-1},v(x^k)))}{\max\left\{\mu\vert h(x^k,d^{k-1}) h(x^{k-1},v(x^k))\vert,-\mu h(x^{k-1},v(x^{k-1}))\vert h(x^{k-1},v(x^k))\vert\right\}}.
\end{equation}
 \item[{\rm Step 5.}] Computing the stepsize $\alpha_k$ by some line search strategies.

\item [{\rm Step 6.}]
 Set $x^{k+1}=x^{k}+\alpha_kd^{k}$ and $k = k+1$, return to Step 2.
\end{enumerate}
\end{algorithm}

\begin{remark}

\begin{enumerate}
\item[(i)] If Algorithm \ref{alg1} stops at iteration $k$, then Lemma \ref{stop} implies that $x^k$ is $K$-Pareto critical point.
\item[(ii)] If $m=1$, $K=\mathbb{R}^+$, and $C=\{1\}$, $\beta_k^{MPRP}$ in \eqref{betak} can be write as
\begin{equation}
\beta_{k}=\frac{g_k^{\top}(\Vert g_{k-1}\Vert g_k-\Vert g_{k}\Vert g_{k-1})}{\max\left\{\mu \Vert g_{k-1}\Vert^3,\mu\Vert g_{k}\Vert\Vert g_{k-1}\Vert\Vert d_{k-1}\Vert\right\}},
\end{equation}
which is the modified version for the conjugate parameter proposed by \cite{28}.
\end{enumerate}
\end{remark}

In this paper, we consider two strategies to find the appropriate stepsize, i.e., the standard Wolfe line search and the Armijo line search. Firstly we state the standard Wolfe line search for vector optimization.
\subsection*{Standard Wolfe line search}
Let $d^k\in\mathbb{R}^n$ be a $K$-descent direction for $F$ at the point $x^k$, and $e\in \text{int} (K)$ a vector such that
\begin{equation}
\langle w,e\rangle\leq 1\quad \text{for all} \quad w\in C.
\end{equation}
It is said that $\alpha_k>0$ satisfies the standard Wolfe line search if
\begin{subequations}
\begin{align}
F(x^k+\alpha_k d^k)&\preceq_K F(x^k)+\rho \alpha_k h(x^k,d^k)e,\label{wolfe1}\\
h(x^k+\alpha_k d^k,d^k)&\geq \sigma h(x^k,d^k),\label{wolfe2}
\end{align}
\end{subequations}
where $0<\rho<\sigma<1$.

And the Armijo line search is defined as follows.
\subsection*{Armijo line search}
Let $d\in\mathbb{R}^n$ be a $K$-descent direction for $F$ at the point $x^k$, $0<\rho<1$, $0<\delta<1$ and $e\in \text{int} (K)$ a vector such that
\begin{equation}\label{armijo-e}
\langle w,e\rangle\leq 1\quad \text{for all} \quad w\in C.
\end{equation}
Set $\tau_k=-\frac{h(x^k,d^k)}{\Vert d^k\Vert^2}$, it is said that $\alpha_k=\max\left\{\tau_k,\delta\tau_k,\delta^2\tau_k,\cdots\right\}$ satisfies the Armijo line search if
\begin{equation}\label{armijo}
F(x^k+\alpha_k d^k)\preceq_K F(x^k)+\rho \alpha_k h(x^k,d^k)e.
\end{equation}
\begin{remark}
\begin{enumerate}
\item[(i)]We claim that the definition of the standard Wolfe line search is the same as the standard Wolfe line search in \cite[definition 3.1]{2}. Actually, since $e\in K$ and $\langle w,e\rangle>0$ is required in the definition of the standard Wolfe line search in \cite{2}, we have $e\in \text{int} (K)$.
\item[(ii)]The definition of the Armijo line search in our paper is slightly different from the Armijo line search in \cite{1}. Moreover, it is a natural vector extension of the Armijo line search for multiobjective optimization in \cite{45}.
\item[(iii)]In multiobjective optimization, where $K=\mathbb{R}_+^m$ and $C$ is the canonical basis of $\mathbb{R}^m$, we usually take $e=\left[1,\dots,1\right]^{\top}\in\mathbb{R}^m$.
\end{enumerate}
\end{remark}

\begin{proposition}\cite[Proposition 3.2]{2}
Assume that $F$ is continuously differentiable, $d$ is a $K$-descent direction for $F$ at $x$, and there exists $\mathcal{B}\in\mathbb{R}^m$,such that
\begin{equation}
F(x+\alpha d)\succeq_K \mathcal{B},
\end{equation}
for all $\alpha>0$. If $C$, the generator of $K$, is finite, then there exist intervals of positive stepsizes satisfying the standard Wolfe line search \eqref{wolfe1}-\eqref{wolfe2}.
\end{proposition}
The Proposition stated above indicates that if $F$ is continuously differentiable and bounded below along the direction $d$, where $d$ is a $K$-descent direction for $F$ at $x$, and $K$ is a finitely generated cone, there exist intervals of stepsizes satisfying the standard Wolfe line search.
\begin{proposition}\cite[Lemma 1]{46}
If $d$ is a $K$-descent direction for $F$ at $x$, then there exist intervals of stepsizes satisfying the Armijo line search \eqref{armijo}.
\end{proposition}
For the convergence analysis to our algorithm, now we display the more stringent condition in connection with the scalar case
\begin{equation}\label{sufdes}
h(x^{k},d^{k})\leq c h(x^{k},v(x^{k})),
\end{equation}
for some $c>0$ and any $k\geq 0$. In vector optimization, we say that a direction $d^{k}\in \mathbb{R}^n$ meets the sufficient descent condition at $x^{k}$ if and only if \eqref{sufdes} holds. Next, we will prove that the search direction generated by Algorithm \ref{alg1} satisfies the sufficient descent condition.
\begin{proposition}\label{sdc}
For arbitrary $k\geq 0$, the search direction $d^k$ is defined by \eqref{dk}, then
\begin{equation}
h(x^k,d^k)\leq(1-\frac{2}{\mu})h(x^k,v(x^k)),
\end{equation}
where $\mu>2$. It means that $d^{k}$ satisfies the sufficient descent condition \eqref{sufdes} at $x^{k}$ with $c=1-\frac{2}{\mu}$.
\end{proposition}
\section{global convergence for the modified PRP-type conjugate method}
As a consequence of Lemma \ref{stop}, Algorithm \ref{alg1} successfully stops if a $K$-Pareto critical point is found. From now on, we assume that the sequence generated by Algorithm \ref{alg1} is infinite. In this section, we will investigate the global convergence of the proposed method. In order to prove the global convergence of the new method, we require the objective function $F(x)$ to satisfy the following assumptions:
\begin{assumption}\label{ass1}
The cone $K$ is generated finitely and there is an open set $\Lambda$ that satisfies $\mathcal{L}:=\{x\in\mathbb{R}^n|F(x)\leq F(x^0)\}\subset \Lambda$, and the Jacobian $JF$ is $L$-Lipschitz continuous on $\Lambda$, i.e.,
\begin{equation}
\Vert JF(x)-JF(y)\Vert \leq L\Vert x-y\Vert.
\end{equation}
\end{assumption}

\begin{assumption}\label{ass2}
If a sequence $\left\{D_k\right\}_{k\in\mathbb{N}}\subset F(\mathcal{L})$ and $D_{k+1}\preceq_K D_k$ for all $k\geq0$, then there is a $D\in\mathbb{R}^m$ such that $D\preceq_K D_k$ for all $k\geq0$, which means that all monotonically nonincreasing sequences in $F(\mathcal{L})$ are bounded from below.
\end{assumption}

Both of the above assumptions are natural extensions of those made for the scalar case, and under Assumption \ref{ass1} and \ref{ass2}, if the stepsize $\alpha_k$ satisfies the standard Wolfe line search, we establish that the iterative form satisfies a condition of Zoutendijk's type, which is important to prove the global convergence of the conjugate gradient method with our parameter $\beta_k$.
\begin{proposition}\label{Z1}
(Zoutendijk's type condition) If Assumptions \ref{ass1} and \ref{ass2} hold, $d^k$ is given by \eqref{dk}, $\{x^k\}$ is generated by Algorithm \ref{alg1} and the stepsize $\alpha_k$ satisfies the standard Wolfe line search, then we have
\begin{equation}
\sum_{k\geq0}\frac{h^2(x^k,d^k)}{\Vert d^k\Vert^2}<\infty.
\end{equation}
\begin{proof}
Considering that $d^k$ given by \eqref{dk} satisfies the sufficient descent condition \eqref{sufdes}, we can directly obtain the result from of \cite[Proposition 3.3]{2}.
\end{proof}
\end{proposition}

Next, we will prove the global convergence of Algorithm \ref{alg1} with standard Wolfe line search by contradiction.
\begin{theorem}\label{limv}
If Assumptions \ref{ass1} and \ref{ass2} hold, $d^k$ is given by \eqref{dk}, $\{x^k\}$ is generated by Algorithm \ref{alg1} and the stepsize $\alpha_k$ satisfies the standard Wolfe line search conditions, then we have
\begin{equation}
\liminf_{k\rightarrow\infty}\Vert v(x^k)\Vert=0.
\end{equation}
\end{theorem}

In the next Lemma \ref{Z2}, under Assumption \ref{ass1} and \ref{ass2}, if the stepsize $\alpha_k$ satisfies the Armijo line search, we also establish that the iterative form satisfies a condition of Zoutendijk's type.
\begin{proposition}\label{Z2}
(Zoutendijk's type condition) If Assumptions \ref{ass1} and \ref{ass2} hold, $d^k$ is given by \eqref{dk}, $\{x^k\}$ is generated by Algorithm \ref{alg1} and the stepsize $\alpha_k$ satisfies the Armijo line search, then we have
\begin{equation}
\sum_{k\geq0}\frac{h^2(x^k,d^k)}{\Vert d^k\Vert^2}<\infty.
\end{equation}
\end{proposition}

Now that the condition of Zoutendijk's type for Algorithm \ref{alg1} with Armijo line search is obtained, it is easy to prove the global convergence of Algorithm \ref{alg1} with Armijo line search in the sense of $\liminf_{k\rightarrow\infty}\Vert v(x^k)\Vert=0$.
\begin{theorem}
If Assumptions \ref{ass1} and \ref{ass2} hold, $d^k$ is given by \eqref{dk}, $\{x^k\}$ is generated by Algorithm \ref{alg1} and the stepsize $\alpha_k$ satisfies the Armijo line search, then we have
\begin{equation}
\liminf_{k\rightarrow\infty}\Vert v(x^k)\Vert=0.
\end{equation}
\begin{proof}
Consider the result of Proposition \ref{Z2}, the proof is the same as that of Theorem \ref{limv}.
\end{proof}
\end{theorem}

\section{Numerical Experiments}
In this section, we present some numerical experiments, in order to illustrate the potential practical advantages of our proposed method. We compare our modified Polak-Ribi\`{e}re-Polyak-type conjugate gradient method using Wolfe conditions (MPRP-W) and modified Polak-Ribi\`{e}re-Polyak-type conjugate gradient method using Armijo condition (MPRP-A) with the PRP conjugate gradient method, PRP+ conjugate gradient method, and FR conjugate gradient method proposed by \cite{2}. All codes are written in double precision Fortran 90. All of the tested problems are classic in the multiobjective optimization literature, and we assume that $e=[1,1,\cdots,1]^{\top}\in \mathbb{R}_+^m$, $C$ is considered as the canonical basis of $\mathbb{R}_+^m$, and $K=\mathbb{R}_+^m$.

The conjugate gradient methods considered in numerical experiments are as follows:
\begin{enumerate}
\item[(i)] PRP conjugate gradient method: It is implemented using strong Wolfe line search conditions, i.e., the stepsize is obtained by finding a $\alpha_k>0$ such that
\begin{subequations}
\begin{align}
F(x^k+\alpha_k d^k)&\preceq_K F(x^k)+\rho \alpha_k h(x^k,d^k)e,\\
\vert h(x^k+\alpha_k d^k,d^k)\vert&\leq \sigma \vert h(x^k,d^k)\vert.
\end{align}
\end{subequations}
And the conjugate parameter is defined by
\begin{equation}\label{prp}
\beta_k^{PRP}=\frac{-h(x^k,v(x^k))+h(x^{k-1},v(x^k))}{-h(x^{k-1},v(x^{k-1}))}.
\end{equation}

\item[(ii)] PRP+ conjugate gradient method: It is implemented using strong Wolfe line search conditions, and the conjugate parameter is defined by $\beta_k^{PRP+}=\max\{0,\beta_k^{PRP}\}$, where $\beta_k^{PRP}$ is given by \eqref{prp}.
\item[(iii)] MPRP-W conjugate gradient method: The stepsize is obtained by standard Wolfe line search conditions \eqref{wolfe1}-\eqref{wolfe2}. The conjugate parameter $\beta_k^{MPRP}$ is defined by \eqref{betak}, and we take $\mu=2.4$ in numerical experiments.
\item[(iv)] MPRP-A conjugate gradient method: The stepsize is obtained by Armijo line search condition \eqref{armijo}. The conjugate parameter $\beta_k^{MPRP}$ is defined by \eqref{betak}, and we take $\mu=2.4$ in numerical experiments.
\item[(v)] FR conjugate gradient method: The stepsize is obtained by strong Wolfe line search conditions, and the parameter is defined by $\beta_k=\frac{h(x^k,v(x^k))}{h(x^{k-1},v(x^{k-1}))}$.
\end{enumerate}

According to Lemma \ref{stop}, we know that $\theta(x)=0$ if and only if $x\in \mathbb{R}^n$ is a $K$-critical Pareto point of $F$,  so we consider the stop condition and claim convergence when $\theta(x^k)\geq-5\times \text{eps}^{1/2}$, where $\text{eps}= 2^{-52}\approx 2.22\times10^{-16}$ and corresponds to the machine precision. Alternatively, the process terminates and claims failure if the maximum number of iterations, 5000, is reached.

To intuitively feel the advantages and disadvantages of the numerical performance of different algorithms, the numerical comparisons will be presented using performance profiles \cite{43}. And for the sake of completeness, we will briefly explain the performance profile here.  Let $S$ be the set of solvers, $p$ be the set of problems,  and $t_{p,s}$ be the performance (for example, we consider the following performance measurement: number of function evaluations, number of gradient evaluations, CPU time and number of iterations) of the solver $s\in S$ on the problem $p\in P$. We emphasize that lower values of $t_{p,s}$ mean better performances.  The performance ratio is $r_{p,s}:=t_{p,s}/\min\{t_{p,s}\mid s\in S\}$, and the cumulative distribution function $\rho_s:[1,\infty)\rightarrow [0,1]$ is
\begin{equation}
\rho_s(\tau)=\frac{\vert\{p\in P\mid r_{p,s}\leq \tau\}\vert}{\vert P\vert}.
\end{equation}
Note that $\rho_s(1)$ means the probability that the solver defeats the remaining solvers and is the most efficient over all the considered algorithms. And we can compare the different methods with respect to robustness rates which are readable on the right vertical axes of the associated performance profiles.

\begin{figure}[htbp]
  \centering
  \subfloat[Performance in terms of Number of iterations]
  {\includegraphics[width=0.5\textwidth]{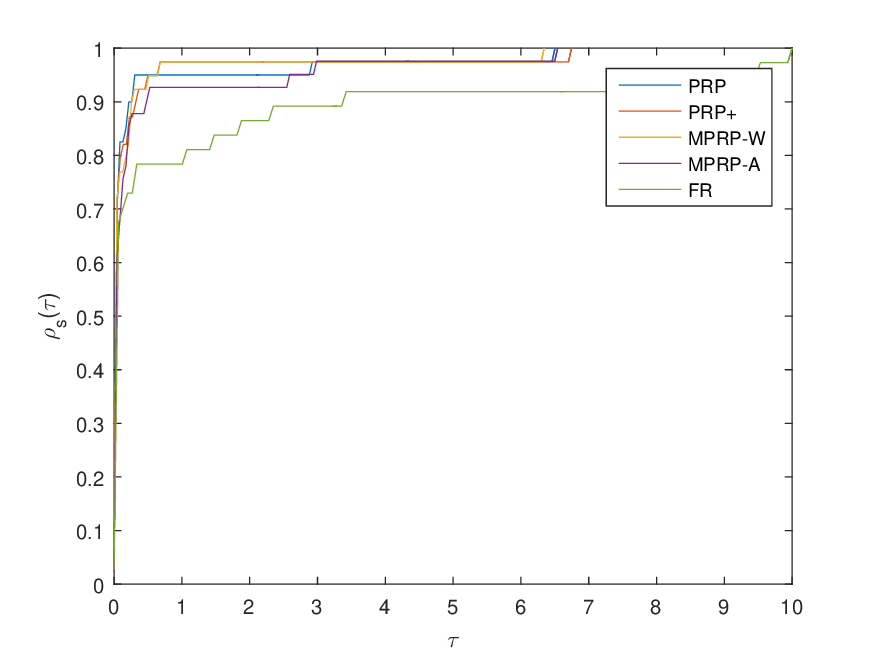}\label{subfig1}}
  \subfloat[Performance in terms of CPU time]
  {\includegraphics[width=0.5\textwidth]{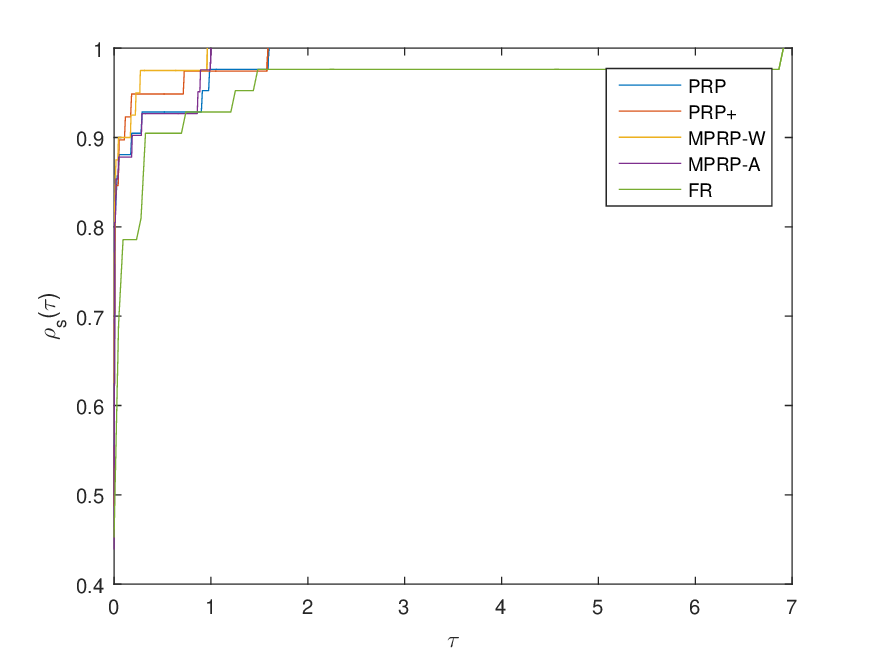}\label{subfig2}}

  \subfloat[Performance in terms of Number of function evaluations]
  {\includegraphics[width=0.5\textwidth]{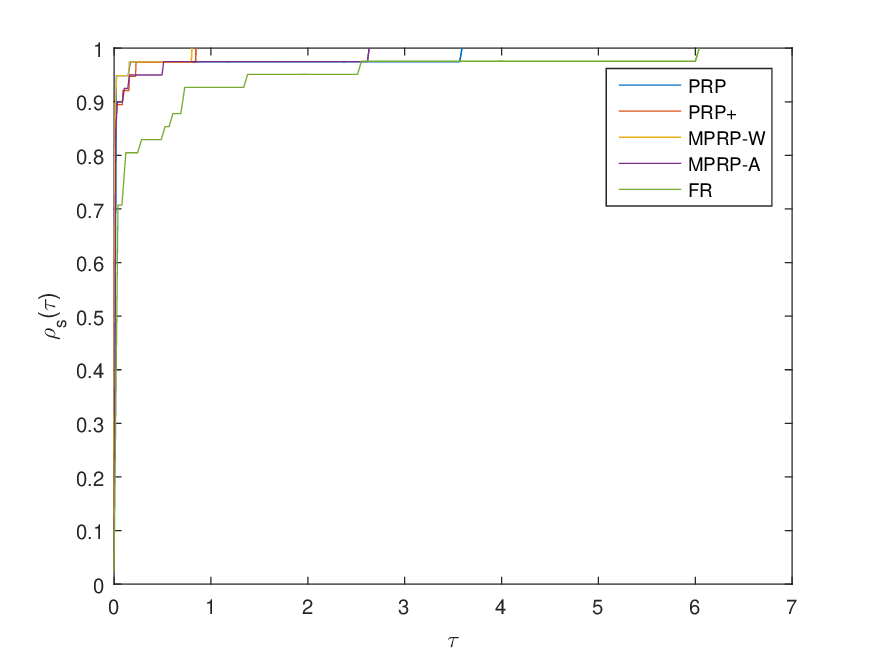}\label{subfig3}}
  \subfloat[Performance in terms of Number of gradient evaluations]
  {\includegraphics[width=0.5\textwidth]{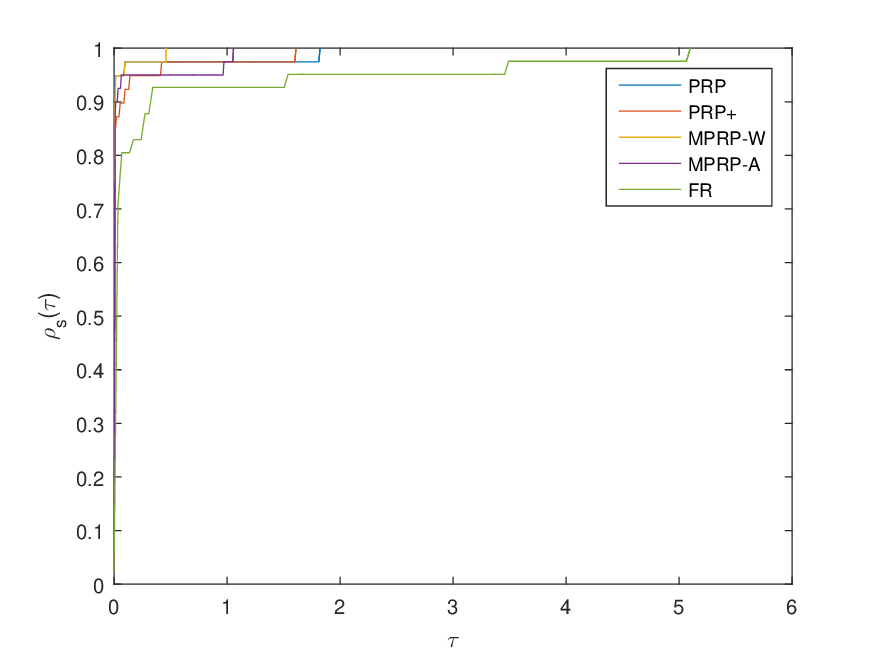}\label{subfig4}}

  \caption{Performance profiles using 200 initial points for each test problem considering the following performance measurement: (A) Number of iterations; (B) CPU time; (C) Number of function evaluations; (D) Number of gradient evaluations.}
\label{fig1}
\end{figure}

As we can see from Figure \ref{fig1}, overall, the four PRP methods are clearly superior to the FR methods in terms of various performance measurement. With respect to the number of iterations (Figure \ref{subfig1}), the MPRP-W method is the most efficient algorithm followed by the PRP method, the MPRP-A method, the PRP+ method, and the FR method. Regarding the CPU time (Figure \ref{subfig2}), although the PRP+ method is slightly the most efficient, it is quickly outperformed by the MPRP-W method, moreover, the MPRP-W method was the first to reach $\rho_s(1)$, which shows that it defeats the remaining method and is the most efficient. In term of the number of function evaluations (Figure \ref{subfig3}), the MPRP-W method is the most efficient and robust one followed by the PRP+ method, the MPRP-A method, the PRP method, and of these comparisons, the FR method performs the worst. Considering the number of gradient evaluations (Figure \ref{subfig4}), the most efficient method is also the MPRP-W method, the MPRP-A method outperforms the PRP+ method, the PRP method and the FR method to be the second most superior method in this measurement. This behaviour is justified by the fact that it generally requires a reasonable number of iterations and the implementation of its backtracking procedure does not use any additional derivative information. Figure \ref{fig1} shows that the MPRP-W method performs well under all performance measurement and is an efficient way to find Pareto points.

From the experiments stated above, it can be seen that the MPRP-W method performs quite well under the four performance measurement: (A) Number of iterations; (B) CPU time; (C) Number of function evaluations; (D) Number of gradient evaluations.

\section{Conclusion}
In this paper, we proposed and analyzed a nonlinear modified Polak-Ribi\`{e}re-Polyak type conjugate gradient method with a nonnegative conjugate parameter to find critical points of vector-valued functions with respect to the partial order induced by a closed, convex, and pointed cone with nonempty interior. This variant are nontrivial extensions of a new Polak-Ribi\`{e}re-Polyak type method of the scalar case to the vector setting. We showed that the search direction in our method satisfies the sufficient descent condition independent of any line search. Furthermore, under mild assumptions, we obtained the results of global convergence with the standard Wolfe line search conditions as well as the standard Armijo line search strategy without convexity assumption of the objective functions. Numerical experiments showed that the effectiveness of the proposed method.

\end{document}